\documentclass[12pt]{article}

\usepackage{amsmath,amssymb}

\newtheorem{theorem}{Theorem}
\newtheorem{lemma}[theorem]{Lemma}
\newtheorem{prop}[theorem]{Proposition}
\newtheorem{corr}[theorem]{Corollary}

\newcommand{\D}{{\mathcal{D}}}
\newcommand{\K}{{\mathcal{K}}}

\newcommand{\C}{{\mathbb{C}}}
\renewcommand{\O}{{\mathcal{O}}}
\renewcommand{\H}{{\mathbb{H}}}

\title{Deformation of line bundles on coisotropic subvarieties.}
\author{Vladimir Baranovsky, Victor Ginzburg, Jeremy Pecharich}

\date{September 29, 2009}

\begin{document}

\maketitle

\begin{abstract}
We prove a criterion stating when a line bundle on a smooth 
coisotropic subvariety $Y$ of a smooth variety $X$ with 
algebraic Poisson structure, admits a first order deformation 
quantization. 
\end{abstract}

\section{Introduction} 

In this paper, partially motivated by \cite{BG}, we consider a smooth
algebraic variety $X$ with an algebraic Poisson structure $P \in H^0(X, \Lambda^2 T_X)$ and a smooth subvariety $Y \subset X$. 
Any line bundle $L$ on $Y$ defines a sheaf of 
$\mathcal{O}_X$-modules. Since $P$ gives a first order deformation 
$\mathcal{A}$ of the structure sheaf $\mathcal{O}_X$, it is natural to ask 
when $L$ can be deformed to an $\mathcal{A}$-module $\mathcal{L}$. In fact,
we consider a slightly more general first order non-commutative deformation
$\mathcal{A}$ which also depends on a class in $H^1(X, T_X)$. 

Based on the standard formalism of hypercohomology and a version
of the deformation complex of Gerstenhaber and Schack, cf. 
\cite{GS} and \cite{MFY},  one expects three
obstruction classes which should vanish if such an $\mathcal{L}$ is to 
exist: these belong to the groups $H^0(Y, \Lambda^2 N)$, $H^1(Y, N)$
and $H^2(Y, \mathcal{O}_Y)$, respectively, where $N$ is the normal
bundle of $Y$ in $X$. The first obstruction is local in nature, and
according to  \textit{loc. cit.} its vanishing simply means that $Y$ 
should  be coisotropic with respect to $P$, i.e. the natural 
projection of $P$ to $\Lambda^2 N$ should be zero. We impose this assumption
on $Y$ throughout this paper and show in Section 3.2 that  in this case
$\mathcal{L}$ always exists Zariski locally.

Next we consider the obstruction class in $H^1(Y, N)$ formulating, 
cf. Theorem 7 of Section 3.3, a 
precise condition on $c_1(L) \in H^1(Y, \Omega^1_Y)$ which guarantees 
existence of a global  deformation $\mathcal{L}$. In general, it will
only be a twisted sheaf of $\mathcal{A}$-modules; one has
an honest sheaf of modules if and only if a further class in 
$H^2(Y, \mathcal{O}_Y)$ vanishes (similarly, coherent sheaves on a general  
algebraic variety $X$ may be twisted by a class in $H^2(X, \mathcal{O}_X)$). 
In Section 4 we assume that $H^2(Y, \mathcal{O}_Y)$ is trivial, but this is
mostly for convenience since one could work with twisted sheaves
of modules instead.

When $X$ is algebraic symplectic, $Y$ is Lagrangian and 
the class in $H^1(X, T_X)$ vanishes, i.e. 
$\mathcal{A} = \mathcal{O}_X \oplus \epsilon \mathcal{O}_X$, 
our condition on $L$ simply says
that $2 c_1(L) = c_1(K_Y)$ in $H^1(Y, \Omega^1_Y)$, cf. Section 5.

If $X$ is symplectic but $Y$ just coisotropic, we can consider the 
standard null foliation $T_F \subset T_Y$ obtained by applying the Poisson 
bivector to the normal bundle of $Y$. In this case, let $L_1$ 
be a line bundle on $Y$ admitting a first order deformation, and $L_2, M$
two line bundles such that $L_2 = M \otimes_{\mathcal{O}_Y} L_1$.
We show that $L_2$ admits a first order deformation if and only if $M$
has a partial algebraic connection along the  null foliation. 
A similar statement is expected for second order deformations if
the partial connection on $M$ is flat (we prove the ``only if" part). 
These results are explained
in Section 4.2 and 4.3, cf. Theorems 10 and 12.

\bigskip
\noindent
\textbf{Acknowledgements.} The first author was partially supported by 
the Sloan Research Fellowship. 

\section{Generalities}

\subsection{First cohomology of a complex.}

Let $\K = \{\K^0 \to \K^1 \to \K^2 \to \ldots\}$ be a complex of sheaves on 
$X$ concentrated in positive degrees. We briefly recall one interpretation
of the hypercohomology group $\H^1(\K)$. 
Let $\mathcal{H}^i = \mathcal{Z}^i/\mathcal{B}^i$
be the cohomology sheaves of $\K$. The standard spectral sequence yields
$$
0 \to H^1(X, \mathcal{H}^0) \to \H^1(X, \K) \to H^0(X, \mathcal{H}^1) 
\to H^2(X, \mathcal{H}^0) \to \ldots
$$
A class in $H^0(X, \mathcal{H}^1)$ is represented by an
open covering $\{U_i\}_{i \in I}$ and  sections 
$\alpha_i \in \Gamma(U_i, \mathcal{Z}^1)$ such that on $U_i \cap U_j$  
$$
\alpha_i - \alpha_j = d \beta_{ij}\qquad \beta_{ij} \in \Gamma(U_i \cap U_j, 
\K^0)
$$
By definition the elements $d\beta_{ij}$ satisfy the cocycle 
condition on triple intersections and
on $U_i \cap U_j \cap U_k$ the expression 
$\gamma_{ijk} = \beta_{ij} + \beta_{jk} + \beta_{ki}$ will be a
section in $\Gamma(U_i \cap U_j \cap U_k, \mathcal{H}^0)$; 
which gives a class in 
$H^2(X, \mathcal{H}^0)$ as in the sequence above. If this class is zero, 
refining $\{U_i\}$ if necessary we can adjust $\beta_{ij}$ 
by adding an element of $\Gamma(U_i \cap U_j, \mathcal{H}^0)$ to ensure 
$$
\beta_{ij} + \beta_{jk} + \beta_{ki} = 0
$$
on $U_i \cap U_j \cap U_k$.
Thus, a class in $\H^1(\K)$ is represented by a covering 
$U_i$, sections $\alpha_i \in \Gamma(U_i, \mathcal{Z}^1)$ and 
$\beta_{ij} \in \Gamma(U_i \cap U_j, \mathcal{K}^0)$ such 
that $\alpha_i - \alpha_j = d\beta_{ij}$ on $U_i \cap U_j$  
and $\beta_{ij} + \beta_{jk} + \beta_{ki} = 0$
on $U_i \cap U_j \cap U_k$. For $\{\beta_i \in \Gamma(U_i, \K^0)\}_{i \in I}$
 the collection 
$\alpha'_i = \alpha_i + d \beta_i$ and $\beta'_{ij} = \beta_{ij} + \beta_i - \beta_j$ represents the same class. 

\subsection{Three examples.}

Let $F_1, \ldots, F_n, G$ be 
sheaves of $\O_X$-modules. We will denote by 
$\D(F_1 \times \ldots \times F_n, G)$ the 
sheaf of algebraic differential operators (of finite order in each of the $n$
variables). It has an increasing filtration by subsheaves 
$\D^k(F_1 \times \ldots \times F_n, G)$
of operators which have total order $\leq k$. 
The $\O_X$-module structure on $G$ gives an $\O_X$-module
structure on differential operators.
If $F_1, \ldots, F_n$ and $G$ are coherent then 
so is $\D^k(F_1 \times \ldots \times F_n, G)$.
We will also write
$F^{\times n}$ for the $n$-fold cartesian product
$F \times \ldots \times F$; and  $\D_0(\O_X^{\times n}, G)$
for the sheaf of differential operators which 
vanish if either of the arguments is a constant.

The three examples described below are sheafifications of a
deformation complex considered by Gerstenhaber and Schack, cf. \cite{GS}. 

\bigskip
\noindent
\textbf{Example 1.}
Set $\K^i(X) = \D_0(\O_X^{\times (i+1)}, \O_X)$ with the
Hochschild differential. Then $\H^1(\K(X)) = H^1(X, T_X) 
\oplus H^0(X, \Lambda^2 T_X)$ parameterizes
flat deformations $\mathcal{A}$ of $\O_X$ over $\C[\epsilon]/\epsilon^2$, 
which locally split as $\O_X \oplus \epsilon \O_X$. Locally these
are given by $\alpha_i \in \Gamma(U_i, \D_0(\O_X \times \O_X, \O_X))$
which satisfy the cocycle condition $d\alpha_i(f, g, h) = \alpha_i(fg, h) -
\alpha_i(f, gh) + \alpha_j(f, g) h - f \alpha_i(g, h) = 0$ and on double
intersections
$$
(\alpha_i - \alpha_j)(fg) = d \beta_{ij}(fg) = \beta_{ij}(fg) - 
f \beta_{ij}(g) - \beta_{ij}(f) g.
$$
In addition, on triple intersections $U_i \cap U_j \cap U_k$ 
we require $\beta_{ij} + \beta_{jk} = \beta_{ik}$.
Write each $\alpha_i(f, g)$ as a sum of its symmetric and
antisymmetric part:
$$
\alpha_i(f, g) = \alpha^+_i(f, g) + \alpha^-_i(f, g)
= \frac{1}{2} (\alpha_i(f, g) + \alpha_i(g, f)) + \frac{1}{2}
(\alpha_i(f, g) - \alpha_i(g, f))
$$
Since $d\beta_{ij}$ is symmetric in $f, g$, we see 
that $\alpha_i^- = \alpha^-_j$
on $U_i \cap U_j$. Moreover, the cocycle condition implies that 
$\alpha_i^-$ is a first order operator in each of the arguments. 
Since $\alpha_i^-$ also vanish on constant functions, they glue into
a section $\alpha^-$ in $H^0(X, \Lambda^2 T_X)$. 

As for the symmetric part, since each $U_i$ is smooth, we can write
$\alpha_i^+ = d\beta_i$ for some $\beta_i \in \Gamma(U_i, \D_0(\O_X, \O_X))$.
Then on every double intersection 
 $\beta_i - \beta_j - \beta_{ij}$ is 
a derivation of $\mathcal{O}_X$. This defines a class in  $H^1(X, T_X)$. 
Conversely, if $\beta_{ij}$ are vector fields on 
$U_i \cap U_j$ representing a class in $H^1(X, T_X)$
then $a_0 + \epsilon a_1 \mapsto a_0 + \epsilon(a_1 + \beta_{ij}(a_0))$
give transition functions which allow to glue the first order deformations
$\O_{U_i} \oplus \epsilon \O_{U_i}$ with cocycle $\alpha^-$, into
a sheaf of algebras $\mathcal{A}$. 

We also observe that $\H^0(\K(X)) = H^0(X, T_X)$ classifies 
those automorphisms of $\mathcal{A}$ which restrict to the identity 
$mod(\epsilon)$.

\bigskip
\noindent
\textbf{Example 2.} If $F$ is a coherent sheaf on $X$ consider
$\K(F)$ with $\K^i(F) = \D_0(\O_X^{\times i}, \D(F, F))$ with the
Hochshild differential corresponding to the natural 
$\O_X$-bimodule structure on $\D(F, F)$, see \cite{We}. 
Then $\H^1(\K(F)) = Ext^1(F, F)$ parameterizies 
flat deformations $\mathcal{F}$ of $F$ to a module over $\O_X[\epsilon]/\epsilon^2$, and $\H^0(\K(F)) = Hom_X(F, F)$ can
be identified with automorphisms of $\mathcal{F}$ which restrict to the
identity $mod(\epsilon)$. 

\bigskip
\noindent
\textbf{Example 3.} The $\O_X$-module structure on $F$ gives 
a morphism of $\O_X$-bimodules 
$\O_X \to \D(F, F)$ hence a morphism of complexes $\mathcal{K}(X) 
\to \mathcal{K}(F)$. Set $\mathcal{K}(X, F) = Cone(\mathcal{K}(X) 
\to \mathcal{K}(F))(-1)$ so that $\mathcal{K}^i(X, f) 
= \D_0(\O_X^{\times (i+1)}, \O_X) \oplus \D_0(\O_X^{\times i}, 
\D(F, F)))$. By the previous subsection and \cite{GS} 
 $\H^1(\K(X, F))$
corresponds to isomorphism classes of flat deformations of $(\O_X, F)$ 
to a pair $(\mathcal{A}, \mathcal{F})$ where $\mathcal{F}$ is a
left $\mathcal{A}$-module. Observe that in \textit{loc. cit.}
the deformation complex has three factors. The
extra factor would correspond in our case 
to deforming the algebra structure of $\D(F, F)$
as well. Since in this paper we will not consider such 
deformations we omit the third factor of the Gerstenhaber-Schack
deformation complex.
 
\bigskip
\noindent
A short exact sequence of complexes $0 \to \K(F) \to \K(X, F) \to \K(X) \to 0$
induces a long exact sequence of cohomology 
\begin{equation}
\label{seq}
\ldots \to \H^1(X, \K(F)) \to \H^1(X, \K(X, F)) \to \H^1(X, \K(X)) 
\to \H^2(X, \K(F)) \to \ldots
\end{equation}
Therefore, given an isomorphism class of $\mathcal{A}$, one expects
that a coherent sheaf $F$ admits a deformation to a left $\mathcal{A}$-module
$\mathcal{F}$ if and only if a certain class $c \in \H^2(X, \K(F))$
vanishes. 

\bigskip
\noindent
The purpose of this note is to make this vanishing condition explicit 
in a particular case. 

\section{Line bundles on coisotropic subvarieties}

From now on we fix a closed embedding $\eta: Y \to X$ of a smooth subvariety,
a first order deformation $\mathcal{A}$ of $\O_X$ with a class
$(\kappa, \frac{1}{2} P) \in H^1(X, T_X) \oplus H^0(X, \Lambda^2 T_X)$, 
and a line bundle $L$ on $Y$. Any vector bundle on $Y$
may also be viewed as $\O_X$-module by applying $\eta_*$ and 
in such a case we abuse notation by dropping $\eta_*$ to make
the formulas more readable. 
We would like to state explicitly 
when $L$ admits a deformation to a left $\mathcal{A}$-module. 
Write $T_Y, N, I$ for the tangent bundle, normal 
bundle and the ideal sheaf of $Y$, 
respectively. 

\subsection{Cohomology}

We observe that the cohomology sheaves of 
$\mathcal{H}^p(\K(X))$ are given by  $\Lambda^{p+1} T_X$ 
due to the Hochschild-Kostant-Rosenberg isomorphism.  
\begin{prop} 
The $p$-th cohomology sheaf $\mathcal{H}^p(\K(L))$ is 
isomorphic to $\Lambda^p N$. 
\end{prop}

\bigskip
\noindent
For $p = 0, 1$ we can make it explicit. First, $\mathcal{H}^0(\K(F)) 
= \mathcal{H}om_{\O_X} (F, F)$ for any sheaf $F$ of 
$\O_X$-modules by definition.
In our case this gives $End_{\O_Y}(L)=\O_Y$. For $p = 1$ 
a section of $\mathcal{H}^1$ is represented locally by a
 map $\alpha_L: \O_X \otimes_{\C} L \to L$ which 
satisfies 
$$
\alpha_L(fg, l) = \alpha_L(f, gl) + f \alpha_L(g, l).
$$
This immediately gives $\alpha_L|_{I^2 \otimes_\C L} = 0$. 
Restricting $\alpha_L$ to $(I/I^2) \otimes_{\C} L$ we see that for 
a local section $x$ on $I$ and a local section $f$ of $\O_X$ we have
$$
f \alpha_L(x, l) = \alpha_L(fx, l) = \alpha_L(xf, l) = \alpha_L(x, fl) 
$$
where the first and the third equalities follow from the cocycle condition. 
Recalling that $I/I^2 = N^\vee$ we can view the restriction of $\alpha_L$
as an $\O_Y$-bilinear map  $N^\vee \otimes_{\C} L \to L$. This gives
 a morphism $\mathcal{H}^1(L) \to  N$. 
To show that it is an isomorphism
it suffices to restrict to an affine open subset on which  $L\simeq \O_Y$
and $Y$ is given by vanishing of a regular sequence;
then the isomorphism follows from the Koszul
complex.
The spectral sequence of hypercohomology for $\K(L)$ gives in lower
degrees
$$
0 \to H^1(Y, \O_Y) \to \H^1(\K(L)) 
\to H^0(Y,  N) \to  H^2(Y, \O_Y) \to \widehat{\H}^2(\K(L)) 
\to H^1(Y, N) \to \ldots
$$
$$
0 \to \widehat{\H}^2(\K(L)) \to \H^2(\K(L)) \to H^0(Y, \Lambda^2 N) \to 
\ldots
$$
where $\widehat{\H}^2(\K(L))$ is a subspace of 
$\H^2(\K(L))$ which may be defined through the second line. 
Therefore existence of a deformation $\mathcal{L}$ of $L$ is equivalent
to the vanishing of a certain class $c \in \H^2(\K(L))$ and the
latter can be split as a chain of conditions:
\begin{equation}
\label{raz}
 \textrm{the image of }c\textrm{ in }
H^0(Y, \Lambda^2 N)
\textrm{ is zero, and}
\end{equation}
\begin{equation}
\label{dva}
\textrm{the image of }c\textrm{ in }
H^1(Y,  N)\textrm{ is zero, and}
\end{equation}
\begin{equation}
\label{tri}
\textrm{the image of }c\textrm{ in }
Coker(H^0(Y,  N)\to H^2(Y, \O_Y)
\textrm{ is zero}
\end{equation}
Note that \eqref{dva} can be formulated only if \eqref{raz} holds, 
and \eqref{tri} can be formulated only \eqref{raz}, \eqref{dva} hold. 
The following lemma follows immediately from the definitions
\begin{lemma}
The image of $(\kappa, \frac{1}{2}P)\in \H^1(\K(X))$ under
the composition 
$$
\H^1(\K(X)) \to \H^2(\K(L))\to H^0(Y, \Lambda^2 N)
$$ 
is equal to the projection of $\frac{1}{2} P \in H^0(X, \Lambda^2 T_X)$
to $H^0(Y, \Lambda^2 N)$. In other words, \eqref{raz} 
holds if and only if $Y$ is coisotropic with respect to the bivector field $P$.
\end{lemma}

\medskip
\noindent
Throughout the rest of the paper we will assume that $Y$ is coisotropic in $X$. 
In this case, there is a well-defined projection of $P$ to 
$H^0(Y, N \otimes_{\mathcal{O}_Y} T_Y)$ which we denote again by $P$.

\subsection{The affine case.}

When $X = Spec(A)$ and $Y = Spec(B)$ are affine the conditions
\eqref{dva} and \eqref{tri}
become trivial. Since $H^1(X, T_X)$ is trivial we can assume
that $\alpha_X(f, g) = \frac{1}{2}P(df, dg)$. 
For $I= Ker(A \to B)$ the $B$-module $N^\vee = I/I^2$ is isomorphic
to the global sections of the conormal bundle to $Y$ in $X$.
In this subsection we write $A$, $B$, $\Omega_A$ and
$\Omega_B$ for the global sections of the sheaves $\mathcal{O}_X$,  
$\mathcal{O}_Y$, $\Omega_X^1|_Y$ and $\Omega^1_Y$, 
respectively. In particular, there is a short exact sequence of $B$-modules
$$
0 \to N^\vee \to \Omega_A\to \Omega_B \to 0
$$

\begin{theorem}
\label{exist-alpha}
For any $Y = Spec(B) \subset X= Spec(A)$ and $L$ as above there 
exists $\alpha_L: A \otimes_\C L \to L$ such that 
\begin{equation}
\label{alpha-el}
\alpha_L(aa', l) - \alpha_L(a, a'l) + \alpha_X (a, a') l - 
a \alpha_L(a', l) = 0
\end{equation}
Any such $\alpha_L$ vanishes on $I^2 \otimes_\C L$.
Moreover, it may be taken in the form
$$
\alpha_L (a, l) = \psi(a) l + \rho(da, l)
$$
where
$$
\psi \in \D^2_0(A, \D^0(L, L)) = \D^2_0(A, B); 
\qquad \rho \in \D^1_0(A, \D^1(L, L)) = Hom_B(\Omega_A, \D^1(L, L))
$$
and both $\psi$ and $\rho$, in addition to vanishing on the constants,
also vanish on $I^2 \subset A$.  
\end{theorem}
\textit{Proof.} The vanishing on $I^2 \otimes_\C L$ is an 
immediate consequence of the equation imposed on $\alpha_L$ and
the coisotropness condition
$\alpha_X(I, I) \subset I$.
Substituting the expression for $\alpha_L$ in terms of $\psi$ 
and $\rho$ and
using $d(a a') = a d(a') + d(a) a'$ we get:
$$
\big(\psi(a a') - \psi(a) a' - a \psi(a')\big) l + \frac{1}{2} P (da, da') l 
+ \big(a' \rho(da, l) - \rho(da, a' l)\big) = 0 
$$
By assumption $\psi(1) = 0$ and $\psi$ has order 2, while $\rho$ has order
1 in the $l$ variable. Therefore, if $\sigma_\psi \in Hom_B (Sym^2 \Omega_A, B)$ and 
$\sigma_{\rho} \in Hom_B (\Omega_A \otimes_B \Omega_B, B)$ are the corresponding principal symbols, then
we must ensure that
$$
\sigma_{\psi} (da \otimes da') + \frac{1}{2}P (da, da') + 
\sigma_{\rho}(da, d (a'|_Y)) = 0
$$
in $B$. Observe that all three terms may be viewed
as $B$-linear homomorphisms 
$$
\Omega_A \otimes_B \Omega_A \to B.
$$ 
The third term
vanishes on $\Omega_A \otimes_B N^\vee$ by its definition.
Hence if $\sigma_\psi$ is known,
existence of $\sigma_\rho$ is equivalent to the condition 
$$
\big(\sigma_\psi + \frac{1}{2} P\big)|_{\Omega_A \otimes_B N^\vee} = 0.
$$
Since $Y$ is coisotropic, $P$ vanishes on $N^\vee \otimes_B N^\vee$
and therefore we should look for $\sigma_\psi$ 
in the submodule $S \subset Hom_B (Sym^2 \Omega_A, B)$
of homomorphisms which vanish on $Sym^2 N^\vee$.  
For this submodule we can write 
a short exact sequence
$$
0 \to Hom_B(Sym^2 \Omega_B, B) \to S \to 
Hom_B(\Omega_B \otimes_B N^\vee, B) \to 0
$$
Since $Y$ is affine and smooth, this sequence splits, and the 
image of $-\frac{1}{2}P$ in its quotient term may be lifted
to some $\sigma_\psi \in S \subset Hom_B(Sym^2 \Omega_A, B)$.
This will ensure the vanishing condition for $\sigma_\psi + \frac{1}{2} P$
and therefore existence of $\sigma_\rho$. 
Finally, since $X$, $Y$ are affine and smooth, symbols can be 
lifted to differential operators. This finishes the proof. 
$\square$

\bigskip
\begin{theorem}
\label{lifts}
Let $\alpha_L$ and $\alpha'_L$ be two maps $(A/I^2) \otimes_\C L \to L$
satisfying the first equation of the previous theorem.
Then $(\alpha_L - \alpha_L')$ vanishes on $(I/I^2) \otimes_\C L$
if and only if there exists $\C$-linear $\beta: L \to L$ such that
$$
\alpha_L(a, l) - \alpha_L'(a, l) = \beta (al) - a \beta (l).
$$
This condition means precisely that the two left $\mathcal{A}$-module 
structures on 
$\mathcal{L} = L \oplus \epsilon L$ defined by $\alpha_L$ and $\alpha_L'$, 
respectively, are equivalent via the isomorphism 
$$
l_1 + \epsilon l_2 \mapsto l_1 + \epsilon(\beta(l_1) + l_2).
$$
Moreover, if both $\alpha_L$ and $\alpha_L'$ are bidifferential operators
as in the previous theorem then $\beta \in \D^2(L, L)$. 
\end{theorem}
First we need the following lemma
\begin{lemma}
\label{hoch}
Let $R: A \otimes_\C L\to L$ be a $\C$-linear map. Then 
$$
R(a, l) = \beta(al) - a \beta(l) 
$$
for some $\beta \in Hom_\C(L, L)$ if and only if $R$ vanishes
on $I \otimes_\C A$ and also satisfies
$$
R(ab, l) - R(a, bl) + a R(b, l) = 0
$$
\end{lemma}
\textit{Proof.} Then ``only if" part is obvious. Suppose
that $R$ vanishes on $I\otimes L$, i.e. descends to a linear map
$B \otimes L \to L$ which we denote again by $R$. Then the equation
imposed on $R$ means that $R$ gives a 1-cocycle in the Hochschild complex
of the $B$-bimodule $End_\C(L)$. By Lemma 9.1.9 of \cite{We} we have
$$
HH^i(B, End_\C(L)) = Ext^i_B(L, L).
$$
Since $L$ is projective over $B$ we have $Ext^1_B(L, L) = 0$.
Therefore, $R$ is a coboundary, which means precisely $R(al) = \beta(al) - a \beta(l)$ for $\beta \in End_\C(L)$.
$\square$

\bigskip
\noindent
\textit{Proof of the theorem.} 
The three-term equation of the previous lemma
obviously holds for $\alpha_L - \alpha_L'$. Hence by previous
lemma a required $\xi$ exists but apriori it may be just a linear map.
However, by our choice of $\alpha_L, \alpha_L'$ the difference
is a bidifferential operator which has order $\leq 1$ in $l$. 
Therefore $\beta(l)$ is a differential operator of order $\leq 2$. $\square$

\bigskip
\noindent
Finally, we would like to identify those maps 
$\gamma: N^\vee \otimes_\C L \to L$
which \textit{locally} extend to $\alpha_L$ 
satisfying \eqref{alpha-el}.
Observe that this equation implies 
\begin{equation}
\label{gamma-one}
\gamma(ax, l) + \alpha_X(a, x) l - a \gamma(x, l) = 0; \qquad
a \in A, x \in N^\vee, l \in L
\end{equation}
\begin{equation}
\label{gamma-two}
\gamma(xa', l) - \gamma(x, a'l)  + \alpha_X(x, a') l = 0; \qquad
x \in N^\vee, a' \in A, l \in L
\end{equation}

\medskip

\begin{prop} 
\label{extend}
Any $\gamma: N^\vee \otimes_\C B \to B$ satisfying \eqref{gamma-one}
and \eqref{gamma-two} extends to $\alpha_L: (A/I^2) \otimes_\C L \to L$
satisfying \eqref{alpha-el}. Moreover, $\alpha_L$ may be taken 
in the form $\psi(a)l + \rho(da, l)$ if and only if $\gamma$ has 
the form $\psi'(x) l + \rho'(x, l)$ where $\psi' \in \D^1(N^\vee, B)$
and $\rho' \in Hom_B(I, \D^1(L, L))$. 
\end{prop}
\textit{Proof.} Observe that the conormal sequence 
$0 \to N^\vee \to \Omega_A \to \Omega_B \to 0$ admits a 
$B$-linear splitting
since its terms are projective $B$-modules. 
Let $p$, resp. $q$, be the projectors $\Omega_A \to \Omega_A$ 
such that their images are
identified with $N^\vee$ and $\Omega_B$, respectively.
$$
\alpha_B(a, l) = \gamma(p(da), l) + t(a, l)
$$
where $t \in \D^1_0(A, \D^1(L, L))$ is such that 
$t(a, bl) - b t(a, l) = \alpha_X(q(da), q(db)) l$. Such $t$ may be 
found e.g. by choosing a connection on $L$, which is possible on
any smooth affine variety.

An easy computation shows that \eqref{gamma-one} and \eqref{gamma-two}
imply \eqref{alpha-el}, and that $\gamma$ of the form 
$\psi'(x)b + \rho'(x, b)$ gets lifted to $\alpha_L$ of the form
$\psi(a) b + \rho(da, b)$. Conversely, if $\alpha_L$ is represented 
in such a form then we can take $\psi'$ and $\rho'$ to be 
the restrictions of $\psi$ and $\rho$ to $N^\vee = I/I^2 \subset A/I^2$, 
respectively. For the orders of these operators, we observe that 
whenever $I$ annihilates an $A$-module $M$, any degree $\leq k$
operator $A \to M$ restricts to a degree $\leq (k-1)$ operator
$I \to M$ by an easy induction involving the definition of degree. 
$\square$

\subsection{The obstruction in $H^1(Y, N)$.}

We consider the general situation when the first 
order deformation $\mathcal{A}$ of $\mathcal{O}_X$ is non 
split, i.e. both classes $P \in H^0(X, \Lambda^2 T_X)$ and 
$\kappa \in H^1(X, T_X)$
are nonzero. Thus, we have an open covering $\{U_i\}$ of $X$ and
$\mathcal{A}|_{U_i} \simeq \O_X \oplus \epsilon \O_X$, while the collection
of vector fields $\beta_{X, ij}$ on $U_i \cap U_j$, 
representing the class $\kappa$, gives transition functions
between the two trivializations of $\mathcal{A}|_{U_i \cap U_j}$:
$f_0 + \epsilon f_1 \mapsto f_0 + \epsilon(f_1 + \beta_{X,ij}(f_0))$. The Leibniz
rule for $\beta_{X,ij}$ ensures that the transition function agrees with the product
$$
(f_0 + \epsilon f_1) * (g_0 + \epsilon g_1)=
f_0 g_0 + \epsilon(f_0 g_1 + f_1 g_0 + \alpha_X(f_0, g_0))
$$
where $\alpha_X(f_0, g_0) = \frac{1}{2} P(df_0, dg_0)$.
We will also assume that $N$ and $L$ are trivial on each $Y\cap U_i$. 

\bigskip
\noindent
In this setup, we would like to find a condition which guarantees existence
of a collection $\{\alpha_{L, i}, \beta_{L, ij}\}$ where 
\begin{itemize}
\item $\alpha_{L,i} \in \Gamma(U_i, \D(\O_X \otimes_\C L, L))$
satisfies the condition 
$$
\alpha_{L, i}(fg, l) - \alpha_{L, i} (f, gl) + \alpha_X(f, g) l 
- f \alpha_{L, i} (g, l)= 0
$$
which means that $\mathcal{L} = (L \oplus \epsilon L)|_{U_i}$ is a module
over $\mathcal{A}|_{U_i}$ with respect to the module structure
$$
(f_0 + \epsilon f_1) * (l_0 + \epsilon l_1) =
f_0 l_0 + \epsilon(f_1 l_0 + f_0 l_1 + \alpha_{L, i} (f_0, l_0))
$$
\item $\beta_{L, ij} \in \Gamma(U_i \cap U_j, \D(L, L))$ satisfy 
$$
\alpha_{L, i} (f, l) - \alpha_{L, j} (f, l) = \beta_{L, ij} (fl) 
- f \beta_{L, ij} (l) - \beta_{X, ij}(f) l  
$$
which means that the transition functions $(l_0 + \epsilon l_1) 
\mapsto (l_0 + \epsilon(l_1 + \beta_{L, ij}(l_0)))$ agree with 
the module structure.
\item $\beta_{L, ij}$ satisfy the cocycle condition on triple
intersections, which guarantees that the modules over $U_i$
may be glued into a left $\mathcal{A}$-module $\mathcal{L}$.
\end{itemize}
\begin{theorem}
\label{class}
Let $at(N) \in H^1(Y, End(N)\otimes \Omega^1_Y)$ be the Atiyah class of $N$.
Then existence of $(\alpha_{L, i}, \beta_{L, ij})$ satisfying the 
first two of the three conditions stated above, is equivalent to the
equation in $H^1(Y, N)$:
$$
\big[- at(N) + 2 Id_{N} \otimes c_1(L) \big] 
\cup P + \overline{\kappa} = 0
$$
where $\overline{\kappa}$  stands is the image of $\kappa \in H^1(X, T_X)$
in $H^1(Y, N)$ and $(\cdot) \cup P$ stands for the Yoneda product of 
a class in $H^1(Y, End(N)\otimes \Omega^1_Y) \simeq Ext^1_Y(N \otimes T_Y, N)$
with the image of $P \in H^0(X, \Lambda^2 T_X)$ in $H^0(Y, N \otimes T_Y)$.
\end{theorem}
\textit{Proof.} By Theorem \ref{exist-alpha} we can always find $\alpha_{L, i}$ satisfying
the first equation. To find $\beta_{L, ij}$ with
\begin{equation}
\label{solve}
\alpha_{L, i} (f, l) - \alpha_{L, j} (f, l) + \beta_{X, ij}(f) l
= \beta_{L, ij} (fl) - f \beta_{L, ij} (l),
\end{equation}
first observe that existence of $\beta_{L, ij}$ does not depend on the choice
of $\alpha_{L, i}$ since any other choice will be given by adding 
$\eta_{i}(f, l)$ such that 
$\eta_{i}(fg, l) - \eta_{i} (f, gl) + f \eta_{i}(g, l) = 0$.
We have assumed that $L$ is identified with $\O_Y$ on $U_i$ and 
hence for any section $l$ of $L$ on $U_i$ we can find a 
$\C$-linear splitting $l \mapsto \widehat{l}$ of the 
surjection $\O_X \to \O_Y \simeq L$ on $U_i$. Then 
$$
\eta_i(f, l) = \eta_i(f, \widehat{l} \cdot 1) = \eta_i(f\widehat{l}, 1) - f \eta_i(\widehat{l}, 1)
$$
tells us that if $\{\beta_{L, ij}\}$ solve \eqref{solve} for $\alpha_{L, i}$ 
then $\beta_{L, ij} + \eta_i(\widehat{l}, 1) - \eta_j(\widehat{l}, 1)$
solve the same equation for $\alpha_{L, i} + \eta_i$. 

\bigskip
\noindent
Denote by 
$R(f, l)$ the right hand side of \eqref{solve}. 
If $R|_{I \otimes L} = 0$ then we can apply Lemma \ref{hoch}  and 
find $\beta_{L, ij}$.  
It remains to establish whether we can replace $\alpha_{L, i}$ 
by $\alpha_{L, i} + \eta_i$ so that the right hand side of 
\eqref{solve} vanishes on $I \otimes_\C L$. In other words, we would like
to have equality of maps $I \otimes_\C L \to L$:
$$
\alpha_{L, j} (f, l) - \alpha_{L, i} (f, l) + \beta_{X, ij}(f) l
= \eta_i(f, l) - \eta_j(f, l)
$$
Since $Y$ is coisotropic, each term on
the left hand side vanishes on $I^2 \otimes_\C L$. 
Therefore recalling $I/I^2 \simeq N^\vee$ we can view
the above equality as equality of functions on $N^\vee \otimes_\C L$. 
We observe that if $\eta_i$ are found as functions on $N^\vee \otimes_\C L$ 
we can always extend them to $(\mathcal{O}_X/I^2) \otimes_\C L$
as in Proposition \ref{extend}.

We will show that the left hand side is an $\mathcal{O}_Y$-bilinear map
$N^\vee \otimes L \to L$ for a particular choice of $\alpha_{L, i}$.
The condition $\eta_{i}(fg, l) - \eta_{i} (f, gl) + f \eta_{i}(g, l) = 0$
implies that $\eta_i$ also must be $\mathcal{O}_Y$-bilinear on 
$N^\vee \otimes_{\C} L$. Hence existence of $\{\eta_i\}$ 
will be equivalent to vanishing of a class in $H^1(Y, N)$. 

\bigskip
\noindent
To calculate the class in $H^1(Y, N)$ explicitly, let $x^i_1, \ldots, x^i_r$
be the basis of $N^\vee|_{U_i}$ and $e^i$ a section spanning $L|_{U_i}$. 
On $N^\vee \otimes_\C L$ we can set $\alpha_{L, i}(x^i_s, e^i) = 0$ which
implies 
$$
\alpha_{L, i} (\sum_s a_s x^i_s, b e^i) = 
(\sum_s \alpha_X(x_s^i, a_s) b + 2 \alpha_X (\sum_s a_s x_s, b)) e^i.
$$ 
On a double intersection we have $x_s^i = \sum_s A^{ij}_{sp} x^j_p$
where $A^{ij}$ is the transition matrix. Similarly $e^i = B^{ij} e^j$.
Rewriting $\alpha_{L, j}$ in the basis $x_s^i$ we find that
$$
(\alpha_{L, j} - \alpha_{L, i})(\sum_s a_s x^i_s, b e^i)
= (\sum_{s, p, r} a_s \alpha_X(x^i_r, A^{ij}_{sp})A^{ji}_{pr}
+ 2 \alpha_X (\sum_s a_s x_s^i, B^{ij})B^{ji})) b e^i
$$
By a similar calculation involving lifts of vector fields to elements
of Atiyah algebras we find that $d A^{ij} \cdot A^{ji}$ represents
minus the Atiyah class of $N$ and $d B^{ij} \cdot B^{ij}$ the first
Chern class of $L$. It is clear that the term $\beta_{X, ij}(f) l$ in 
\eqref{solve} represents the class $\overline{\kappa}$ 
as in the statement of the theorem. This finishes the proof.
$\square$

\bigskip
\noindent
\textbf{Remark.} Even if the class in $H^1(Y, N)$ vanishes and
$\beta_{L, ij}$ exist, they may not satisfy the cocycle condition on
triple intersections. However, equation \eqref{solve} implies that
on $U_i \cap U_j \cap U_k$ the expression $\beta_{L, ij} + \beta_{L, jk}
+ \beta_{L, ki}$ is $\mathcal{O}_Y$-linear and thus defines a class 
in $H^2(Y, \mathcal{O}_Y)$. The vanishing of this class,
or a weaker condition \eqref{tri}, is needed to 
ensure that $\mathcal{L}$ exists.

\begin{corr}
Let $Y$ be a coisotropic smooth subvariety in $X$ with 
$H^2(Y, \O_Y) = 0$, $\mathcal{A}$ a first
order deformation of $\O_X$ with class $(\kappa, P)$, and 
$L$ a line bundle on $Y$ such that
$$
\big[- at(N) + 2 Id_{N} \otimes c_1(L) \big] 
\cup P + \overline{\kappa} = 0
$$
in $H^1(Y, N)$. Then $L$ admits a first order deformation $\mathcal{L}$
to a left $\mathcal{A}$-module. If $H^1(Y, \mathcal{O}_Y) = 0$
the set of isomorphism classes 
of such $\mathcal{L}$ is parameterized by $H^0(Y, N)$. In 
general, the group of
automorphisms (restricting to the identity $mod(\epsilon)$)
of each $\mathcal{L}$ is isomorphic to $H^0(Y, O_Y)$. 
\end{corr}
\textit{Proof.} For the isomorphism classes we recall the sequence
\eqref{seq}. By section 3.1
$$
\H^0(\K(L)) = H^0(Y, \O_Y), \quad \H^1(\K(L)) = H^0(Y, N)
$$
Recall that $\H^0(\K(X)) = H^0(X, T_X)$. It follows from the definitions
that $\H^1(\K(X, L))$ is the vector space of all pairs $\partial_L, \partial$, where $\partial_L\in \D^1(L, L)$ and $\partial$
is an extension of the symbol of $\partial_L$ to a vector field on $X$.
It follows that 
$$
0 \to \H^0(\K(L)) \to \H^0(\K(X, L)) \to \H^0(\K(X)) \to 0
$$
is exact. Therefore
$$
0 \to \H^1(\K(L)) \to \H^1(\K(X, L)) \to \H^1(\K(X)) \to \ldots
$$
is also exact, and a lift of any element in $\H^1(\K(X))$
is well defined up to an element of $\H^1(\K(L)) = H^0(Y, N)$. 
More explicitly, a section of $H^0(Y, N)$ restricted to $U_i$
may be lifted to a derivation $\delta_i: \mathcal{O}_X \to \mathcal{O}_Y$
and we can adjust $\alpha_{L, i}(a,l)$ replacing it by 
$\alpha_{L, i}(a, l) + \delta_i (a) l$. On a double
intersection the difference $\zeta_{ij} = \delta_i - \delta_j$
is a derivation $\O_Y \to \O_Y$ which we can view as an operator from $L$
to itself, since $L$ is trivialized on $U_i \cap U_j$. 
Thus the data $(\alpha_{L, i}, \beta_{L, ij})$ will be 
replaced by the data $(\alpha_{L, i} + \delta_i \cdot Id_L, 
\beta_{L, ij} + \zeta_{ij})$

To prove the assertion about automorphisms of $\mathcal{L}$: let
$x$ be such automorphism, then $Id_\mathcal{L} - x$ is an endomorphism
of $\mathcal{L}$ which takes values on $\mathcal{L}/\epsilon \mathcal{L}
\simeq L$ and vanishes on $\epsilon\mathcal{L} \simeq L$, i.e.
a morphism of sheaves $L \to L$. It is easy to see that such morphism
must be $\O_X$-linear, i.e. given by an element of $H^0(Y, \O_Y)$. 
$\square$

\bigskip
\noindent
\textbf{Remark} It follows from the definitions that for 
a deformation $\mathcal{A}$ constructed from a pair $(\kappa, \frac{1}{2}P)$
the deformation $\mathcal{A}^{op}$  corresponds to the  pair
$(\kappa, -\frac{1}{2}P)$. Therefore  a line
bundle $L$ admits a deformation to a \textit{right} 
$\mathcal{A}$-module precisely when
$$
\big[- at(N) + 2 Id_{N} \otimes c_1(L) \big] 
\cup P - \overline{\kappa} = 0.
$$
If $\overline{\kappa} = 0$ then both deformations exist. However
even in the split affine case when $\mathcal{L} = L\oplus \epsilon L$ 
the left and the right  structures can be made to commute
(after appropriate choices of maps $\O_X \times L\to L$ and
$L \times \O_X \to L$)
only if the image of $P$ in $\Lambda^2 T_X|_Y$ belongs
to the subspace $\Lambda^2 T_Y$, i.e. the bivector $P$ should be
rather degenerate along $Y$.

\section{The case $\kappa =0$.}

In this section we assume that $\kappa = 0$, i.e. 
there exists a global splitting
$\mathcal{A} = \O_X \oplus \epsilon \O_X$, 
and that $H^1(Y, \O_Y) = H^2(Y, \O_Y) = 0$. The 
last condition automatically ensures
\eqref{tri}.

\bigskip
\noindent
\textbf{Remark.} With some minor modifications, the arguments below can
be adjusted to the slightly more 
general case when $\overline{\kappa}=0$. This means that 
vector fields $\beta_{X, ij}$ representing the class $\kappa$, 
may be chosen to satisfy $\beta_{X, ij}(I)\subset I$. We leave the 
details to the motivated reader. 

\subsection{Equivalence classes via a global operator.}

Assume that $L$ admits a (non-split) first order deformation
$\mathcal{L}$ to a left $\mathcal{A}$-module. Embedding 
$I \subset \O_X \subset \mathcal{A}$ we get a globally defined map 
$\gamma: I \otimes_\C L \to L$ given by 
$$
x * l = 0 + \epsilon \gamma(x, l)
$$
Repeating the reasoning of Section 3.2 we see that
$\gamma$ descends to $(I/I^2) \otimes_\C L \simeq N^\vee \otimes_\C L$
and satisfies \eqref{gamma-one} and \eqref{gamma-two}. However, since 
$$
\D^1(N^\vee, \O_Y) \cap Hom_{\O_Y}(N^\vee, \D^1(L, L)) = Hom_{\O_Y} (I/I^2, \O_Y) \simeq N,
$$
the splitting $\psi(x) l + \rho(x, l)$ will in general exist only
on the open sets $U_i$ but not globally.
Thus, we
can only say that $\gamma_i$ glue into a global section 
$$
\gamma \in \Gamma(Y, \D^1(N^\vee \times L, L)),
$$
i.e. $\gamma$ has the total order $\leq 1$ in its two arguments.

\begin{prop}
\label{gamma-def1}
Suppose a line bundle $L$ on $Y$ satisfies the condition on $c_1(L)$
stated in Theorem \ref{class}. In the assumption of this section, 
the set of equivalence classes of $\mathcal{A}$-modules $\mathcal{L}$
deforming $L$ is in bijective correspondence with the set of 
globally defined differential operators 
$$
\gamma \in \Gamma(Y, \D^1(N^\vee \times L, L)),
$$
satisfying \eqref{gamma-one} and \eqref{gamma-two}.
\end{prop}
\textit{Proof.} We have seen above that any $\mathcal{L}$ leads
to $\gamma(x, l)$ as in the statement of the theorem. Conversely,
suppose that $\gamma(x, l)$ exists. Taking an affine open covering $\{U_i\}$
 and using the first formula in the proof 
of Proposition \ref{extend} we can extend $\gamma|_{U_i}$ to 
an operator $\alpha_L(a, l): \mathcal{O}_X \otimes L\to L$
defined on $U_i$ and satisfying \eqref{alpha-el}. 
On double intersections $U_i \cap U_j$
the two operators $\alpha_{L, i}$ and $\alpha_{L, j}$ both
extend $\gamma|_{U_i \cap U_j}$ thus their difference 
satisfies the condition of Theorem \ref{lifts} and we can 
find appropriate transition functions $\beta_{L, ij}$ which
automatically satisfy the cocycle condition on triple intersections
due to the assumption $H^2(Y, \O_Y) = 0$. This shows that
the map from equivalence classes to the set of $\gamma$ is onto.

To show that this map is also injective, assume that two deformations
$\mathcal{L}$ and $\widehat{\mathcal{L}}$ are given. We can 
choose a common refinement of the open coverings on which these
deformations split, and assume that they are given by the data
$\{\alpha_{L, i}, \beta_{L, ij}\}$ and $\{\widehat{\alpha}_{L, i}, 
\widehat{\beta}_{L, ij}\}$, respectively. By assumption, on each
$U_i$ both $\alpha_{L, i}$ and $\widehat{\alpha}_{L, i}$ extend 
$\gamma(x, l)$ on $U_i$ and invoking Theorem \ref{lifts} again we can find 
$\beta_{L, i} \in \Gamma(U_i, \D^2(L, L))$ which allows to 
change the splitting of $\widehat{\mathcal{L}}$ in such a way
that $\alpha_{L, i} = \widehat{\alpha}_{L, j}$. Then on double intersections
both $\beta_{L, ij}$ and $\widehat{\beta}_{L, ij}$ solve the
equation $(\alpha_{L, i} - \alpha_{L, j})(f, l) = \beta_{L, ij}(fl) 
- f \beta_{L, ij}(l)$ hence the difference $\beta_{L, ij}$ and $\widehat{\beta}_{L, ij}$ is $\O_Y$-linear, i.e. given by
multiplication of $l$ by a regular function $\widetilde{\beta}_{L, ij}$.
By definition such functions satisfy the cocycle condition on triple
intersections. By our assumption $H^1(Y, \O_Y) = 0$ and we 
can find $\widetilde{\beta}_{L, i} \in \Gamma(U_i, \O_Y)$ such 
that $\widetilde{\beta}_{L, ij}= \widetilde{\beta}_{L, i} 
- \widetilde{\beta}_{L, j}$ on $U_i \cap U_j$. Using $\widetilde{\beta}_{L, i}$
to adjust the splitting of $\widehat{\mathcal{L}}|_{U_i}$ one more
time, we achieve $\beta_{L, ij} = \widehat{\beta}_{L, ij}$. 
This means that the deformations $\mathcal{L}$ and $\widehat{\mathcal{L}}$
are equivalent. $\square$.

\subsection{The case of non-degenerate bivector.}

In this subsection we assume that the bivector $P$ is non-degenerate, 
i.e. gives an isomorphism $\Omega_X \to T_X$. If in addition the
Schouten-Nijenhuis bracket $\{P, P\}$ vanishes, this means that
$X$ has algebraic symplectic structure (but we only need this condition
when discussing the second order deformations).

By non-degeneracy the restriction of $P$ to $Y$ 
embeds $N^\vee$ as a subbundle into $T_Y$. We will denote the
image by $T_F$ and call it the \textit{null foliation} subbundle.
Coisotropness of $Y$ means that $T_F$ is 
\textit{involutive}, i.e. a sheaf of 
Lie subalgebras with respect to the bracket of vector fields.
We define
the \textit{null foliation Atiyah algebra} $At_n(L)$
to be the preimage of $T_F \subset T_Y$ in $\D^1(L, L)$
with respect to the symbol map $\sigma: \D^1(L, L)\to T_Y$, i.e.
first order operators from $L$ to itself with symbol 
in $T_F$. 
Thus we have an extension
$$
0 \to \O_Y \to At_n(L) \to T_F \to 0
$$
Note that $At_n(L)$ is a sheaf of Lie algebras with 
respect to the commutator of differential operators.
Applying the isomorphism $N^\vee \simeq T_F$ we can view
$\gamma$ as a map $T_F \to \D^1(L, L)$
 and equations
\eqref{gamma-one} and \eqref{gamma-two}  - \eqref{gamma-one} 
become
\begin{equation}
\label{gamma-one-null}
\gamma(ax, l) - a \gamma(x, l) = \frac{1}{2} x(a) l
\end{equation}
\begin{equation}
\label{gamma-two-null}
\gamma(x, al) - a \gamma(x, l) = x(a) l.
\end{equation}
The second equation simply says that $\sigma \circ \gamma = Id_{T_F}$. 
However, $\gamma$ is not $\O_Y$-linear, as can be seen from the 
first equation. The meaning of the first equation can be seen 
from the following theorem 
\begin{theorem}
\label{split-class}
For a non-degenerate $P$ the following conditions on the
line bundle $L$ are equivalent
\begin{enumerate}
\item The equation 
$$
\big[- at(N) + 2 Id_{N} \otimes c_1(L) \big] \cup P = 0;
$$ 
holds in $H^1(Y, N)$;
\item There exists $\C$-linear 
splitting $\gamma: T_F \to At_n(L)$ 
satisfying \eqref{gamma-one-null} and \eqref{gamma-two-null}. 
\item There exists an anti-involution $\partial \mapsto 
\partial^*$ on $At_n(L)$ such that 
$$
\sigma(\partial^*) = - \sigma(\partial), \quad 
f^* = f, \quad (f \partial)^* = \partial^* f
$$
where $f \in \O_Y \subset At_n(L)$.
\end{enumerate}
If  $H^2(Y, \mathcal{O}_Y) = 0$ then either of these 
conditions is equivalent to existence
of a first order deformation $\mathcal{L}$ of $L$. 
If in addition $H^1(Y, \O_Y) = 0$ 
the equivalence class of $\mathcal{L}$ is uniquely determined
by $\gamma$.
\end{theorem}
\textit{Proof.} 
Since \eqref{gamma-one-null} and \eqref{gamma-two-null}
are simply reformulations of \eqref{gamma-one} and
\eqref{gamma-two} the equivalence $1 \Leftrightarrow 2$
is essentially proved in Proposition \ref{gamma-def1} 
(note that the vanishing of a class in $H^2(Y, \O_Y)$
is irrelevant to this equivalence).

To prove equivalence $2 \Leftrightarrow 3$ first assume that
$\gamma$ exists and define the $*$-involution to be $+1$
on $\O_Y \subset At_n(L)$ and $-1$ on the image of $\gamma$. 
Conversely, if the $*$-involution exists then its $(-1)$-eigensheaf
projects isomorphically onto $T_F$ and there is a unique
$\gamma$ such that $\gamma \sigma$ is the projection on the
$(-1)$-eigensheaf. A direct easy computation shows that
the conditions imposed on $\gamma$ and $*$ are equivalent.
$\square$

\bigskip
\noindent
Suppose that the first order deformation $\mathcal{A} = \O_X 
\oplus \epsilon \O_X$ extends to a second order deformation 
$\mathcal{A}'$ over $\C[\epsilon]/\epsilon^3$.  Then on 
affine open subsets the product in $\mathcal{A}'$ will be 
given by  
$$
f * g = fg + \epsilon \alpha_X(f, g) + \epsilon^2 \alpha_X'(f, g)
$$
with the usual associativity condition
$$
\alpha_X(a, \alpha_X(b, c)) - 
\alpha_X(\alpha_X(a, b), c) = d \alpha'_X(a, b, c)
$$
By the explicit formula of Section 1.4.2 of \cite{Ko} 
the locally
defined operator $\alpha'_X(f, g)$ may be taken
\textit{symmetric} (after a local choice of an algebraic connection on
the tangent bundle).

\begin{prop}
\label{quasi-flatness}
In the notation the previous theorem, assume that 
$L$ admits a second order deformation to a left
$\mathcal{A}'$-module $\mathcal{L}'$. Then the operator 
$\gamma$ agrees with the Lie brackets:  
$$
[\gamma(\partial_1), \gamma(\partial_2)] 
- \gamma([\partial_1, \partial_2]) = 0
$$
\end{prop}
\textit{Proof.} 
Locally the second order deformation $\mathcal{L}'$ is given by 
$$
a * l = al + \epsilon \alpha_L(a, l)+ \epsilon^2 \alpha_L'(a, l)
$$
for some bidifferential operator $\alpha_L': \O_X \times L\to L$. The usual
associativity equation reads
\begin{equation}
\label{second-el}
\alpha_L(\alpha_X(a_1, a_2), l) - \alpha_L(a_1, \alpha_L(a_2, l)) 
+ \alpha'_X(a_1, a_2) l = \alpha_L'(a_1, a_2l) + a_1 \alpha'_L(a_2, l) - 
\alpha'_L(a_1 a_2, l)
\end{equation}
If $a_1, a_2$ are sections in $I$ then the first two terms on 
the right disappear. Antisymmetrizing in $a_1$ and $a_2$ and using the 
fact that $2\alpha_X: I \times I \to I$ descends to the bracket
of vector fields
on $T_F \simeq I/I^2$, we obtain the result. $\square$

\subsection{Deformations and connections.}

We keep the previous assumptions running: $\kappa = 0$, 
$H^1(Y, \O_Y) = H^2(Y, \O_Y) = 0$ and $P$ is non-degenerate. 

\begin{theorem}
Let $L_1, L_2$ be two line bundles on $Y$ admitting first
order deformations $\mathcal{L}_1$, $\mathcal{L}_2$ 
corresponding to the operators $\gamma_1: T_F \otimes L_1 \to L_1$
and $\gamma_2: T_F \otimes L_2 \to L_2$. Then 
the bundle $M = \mathcal{H}om_{\mathcal{O}_Y} (L_1, L_2)$
admits a partial algebraic connection $\gamma_M: T_F \otimes M \to M$ defined by 
$$
\gamma_M(\partial, \phi)(l_1) = \gamma_2(\partial, \phi_2(l_1))
- \phi(\gamma_1(\partial, l_1))
$$
Conversely, if $M$ is a bundle with a partial connection
$\gamma: T_F \times M \to M$ and $L_1$ is a line bundle admitting
a first order deformation $\mathcal{L}_1$ then the line bundle
$L_2 = M\otimes_{\O_Y} L_1$ also admits a first order deformation 
$\mathcal{L}_2$ with $\gamma_2: T_F \times L_2 \to L_2$
defined by the formula
$$
\gamma_2(\partial, m \otimes l_1) = \gamma_M(\partial, m) \otimes l_1
+ m \otimes \gamma_1(\partial, l_1).
$$
In addition, suppose that $\mathcal{A}$ extends to a second order
deformation $\mathcal{A}'$. If $\mathcal{L}_1$, $\mathcal{L}_2$ 
admit second order deformations to left $\mathcal{A}'$-modules
$\mathcal{L}_1'$, $\mathcal{L}_2'$ then
 the partial connection $\gamma_M$ is flat, i.e.
its curvature in $H^0(Y, \Lambda^2 T^\vee_F)$ is zero. 
\end{theorem}
\textit{Proof.} One needs to show that the first 
formula indeed defines an $\O_Y$-linear operator and 
that the second formula is well-defined on the tensor product over
$\O_Y$, i.e. $\gamma_2(m \otimes al) = \gamma_2(m a \otimes l)$. These
assertions and the assertions about first order deformations 
 follow from \eqref{gamma-one-null} and
\eqref{gamma-two-null} by a straightforward computation.

For the second order deformation, 
let $\partial_1$, $\partial_2$ be two sections of $T_F$. 
We need to show that
$$
\gamma_M(\partial_1, \gamma_M(\partial_2, \phi))
- \gamma_M(\partial_2, \gamma_M(\partial_1, \phi)) 
- \gamma_M([\partial_1, \partial_2], \phi) = 0
$$
for any section $\phi$ of $M = \mathcal{H}om_{\O_Y} (L_1, L_2)$. 
This follows immediately from the formula of previous proposition
and definition of $\gamma_M$.
$\square$

\bigskip
\noindent
\textbf{Remark.} We also expect that, conversely,  if $\gamma_M$ 
is a flat algebraic connection along the null foliation and
$L_1$ admits a second order deformation, then $L_2 = M \otimes L_1$
also admits a second order deformation. In fancier terms, the 
category of second order deformations of line bundles on $Y$ should be 
a \textit{gerbe} over the Picard category of line bundles 
with a flat algebraic connection along the null-foliation.

\section{The Lagrangian case.}

In this section we assume that the bivector $P \in H^0(X, \Lambda^2 T_X)$
is non-degenerate everywhere on $X$ and that $Y$ is Lagrangian, i.e.
its dimension is half the dimension of $X$. Since 
at this moment we work with first order deformations, we will not need the 
condition that the algebraic 2-form defined by $P$ is closed.
The restriction of $P$ to $Y$ defines an isomorphism 
$N \simeq \Omega^1_Y$. Therefore we can write 
$$
at(N)\in H^1(Y, End(N)\otimes \Omega^1_Y) = H^1(Y, \Omega^1 \otimes T_Y
\otimes \Omega^1_Y)
$$
where we agree that the first two factors in the last expression represent
$End(\Omega^1_Y)$. 

\begin{corr}
If $Y \subset X$ is Lagrangian, $H^2(Y, \O_Y) = 0$
and $L$ is a line bundle on $Y$ then 
$L$ admits a first order deformation to a left $\mathcal{A}$-module
is and only if
$$
- c_1(K_Y)+ 2 c_1(L) + \overline{\kappa} = 0
$$
in $H^1(Y, \Omega^1)$. In particular, if $\kappa = 0$ the deformation 
exists if and only if $2c_1(L) = c_1(K_Y)$.
\end{corr}
\textit{Proof.} Since in the Lagrangian case $P \in H^0(Y, N \otimes T_Y)
\simeq H^0(Y, \Omega^1_Y \otimes T_Y)$ is simply the canonical identity
element, the cup product $\big[- at(N) + 2 Id_{N} \otimes c_1(L) \big] 
\cup P$ of the Theorem \eqref{class} simply amounts to the contraction of
a class in $H^1(Y, \Omega^1 \otimes T_Y \otimes \Omega^1_Y)$ 
in the \textit{last two} factors. For $2 Id_{N} \otimes c_1(L)$
this immediately gives $2 c_1(L)$. For $at(N)$ this would give
$c_1(N) = c_1(\Omega^1_Y) = c_1(K_Y)$ by one of the
equivalent definitions of $c_1$ if we were 
contracting in the \textit{first two} factors. However, by 
Proposition 2.1.1 in \cite{Ka} the Atiyah class $at(N) = at(\Omega^1_Y) 
= - at(T_Y)$ is \textit{symmetric} with respect to the first and the
third factors of $\Omega^1 \otimes T_Y \otimes \Omega^1_Y$. 
This finishes the proof. $\square$

\bigskip
\noindent
We also remark that in the Lagrangian case the null foliation bundle $T_F$
is equal to the full tangent bundle $T_Y$, and the partial connections
considered in Section 4.3 become connections in the usual sense.

\bigskip
\small{

\noindent
\textsl{Addresses:\\
VB and JP: 340 Rowland Hall,
University of California - Irvine, Irvine, CA 92697, USA}\\
\textsl{VG:   
5734 S. University Avenue, University of Chicago, Chicago, IL 60637, USA}\\

\noindent
\textsl{Email: \\
vbaranov@math.uci.edu, ginzburg@math.uchicago.edu, jpechari@math.uci.edu}}

\end{document}